\documentstyle{article}
\input amssym.def
\input amssym
\newtheorem{th}[subsection]{Theorem}
\newtheorem{lem}[subsection]{Lemma}
\newtheorem{prop}[subsection]{Proposition}
\newtheorem{cor}[subsection]{Corollary}
\newtheorem{exam}[subsection]{Example}
\newcommand{\F}{{\cal F}}
\newcommand{\C}{{\cal C}}
\newcommand{\M}{{\cal M}}
\newcommand{\G}{{\cal G}}
\newcommand{\cH}{{\cal H}}
\newcommand{\V}{{\cal V}}
\begin{document}
\author{A.A.Davydov}
\title{Pentagon equation and matrix bialgebras}
\maketitle
\begin{center}
{\it Dedicated to A.N.Tyurin on the occasion of his 60-th birthday}
\end{center}
\begin{center}
Department of Mathematics, Moscow State University,\\
 Moscow 119899, Russia
\end{center}
\begin{center}
e-mail:davydov@mech.math.msu.su
\end{center}
\begin{abstract}
We classify coproducts on matrix algebra in terms of solutions to some modification of 
pentagon equation. The construction of Baaj and Skandalis describing finite dimensional 
unitary solutions of pentagon equation is extended to the non-unitary case. We establish 
the relation between Hopf-Galois algebras and solutions to modified pentagon equation.
\end{abstract}
\tableofcontents
\section{Introduction}

This is the first article in the series of papers on algebraic constructions related to 
$6j$-symbols for finite fusion rules. The subject of this article can be considered as the 
case of one element fusion rule. Instead of starting with categorical setting (which is 
treated in the last section) we consider equivalent problem of describing bialgebra structures 
on matrix algebra and reformulate it as a problem of solving some equation closely connected 
with the so-called pentagon equation (see \cite{sk}). The set of solutions to this modified 
pentagon equation projects onto the set of solutions to ordinary pentagon equation. We 
extend Baaj and Skandalis' description of unitary solutions to pentagon equation \cite{sask} 
to the general case which says that any solution corresponds up to a multiplicity to some 
Hopf algebra. We give an algebraic description of solutions to modified pentagon equation. 
In particular, we establish the relation between Hopf-Galois algebras and solutions to 
modified pentagon equation when the corresponding solution to pentagon equation has 
multiplicity one.

The work was partially supported by RFBR grant no. 99-01-01144.

\section*{Acknowledgment}

The work was started during my visit of Macquarie University (Sydney, Australia) and was completed in Max-Planck-Institut f\"ur Mathematik (Bonn, Germany). I would like to thank these institutions for hospitality and inspiring atmosphere. Special thanks to Prof. Ross Street who pointed out to me the references \cite{sk,sask} and explained his own work \cite{st}. 

\section{Matrix bialgebras and pentagon equation}

In this section we give a description of bialgebra structures on the endomorphism algebra $End(V)$ of 
vector space $V$ in terms of some tensors. 
Remind that a {\em bialgebra} is an algebra $E$ together with a unital (identity preserving) homomorphism of algebras $\Delta :E\to E\otimes E$ ({\em comultiplication}) satisfying the so-called {\em coassociativity} axiom: $(\Delta\otimes I)\Delta = (I\otimes\Delta)\Delta$. 
We start with the following auxiliary statement.
\begin{lem}\label{hea}
Let $f:End(V)\to End(U)$ be a unital homomorphism of endomorphism algebras. Then there exists 
a vector space $W$ and an isomorphism $U \simeq V\otimes W$. The composition $\tilde f:End(V)\to End(V\otimes W)$ of 
$f$ with the factorization isomorphism has a form $\tilde f(x) = x\otimes I$.

The centralizer of the subalgebra $End(V)\otimes I$ in the algebra $End(V\otimes W)$ 
coincides with $1\otimes End(W)$. 
\end{lem}
Proof:

The algebra $End(V)$ is simple. Hence
any its module is a direct sum of simple one. 
The natural 
$End(V)$-module $V$ is simple. The homomorphism $f$ defines $End(V)$-module structure on vector space $U$. 
Now evaluation 
map $V\otimes Hom_{End(V)}(V,U)\to U$ gives an isomorphism. Moreover this isomorphism is an 
isomorphism of $End(V)$-modules and the action of $End(V)$ on $V\otimes Hom_{End(V)}(V,U)$ 
has the form $x(v\otimes l) = x(v)\otimes l$. Denoting $Hom_{End(V)}(V,U)$ by $W$ we get the first part of the statement.

The second part of the lemma follows from the identifications:
$$C_{End(V\otimes W)}(End(V)\otimes I) \cong End_{End(V)}(V\otimes W) = End(W).$$
$\Box$

For a linear map $F:V\otimes M\to V\otimes V$ define the maps $F_{12},F_{13},F_{23}$ which are map $F$ acting on first and second (first and third, second and third) tensor functor of appropriate tensor product. 
\begin{th}\label{cea}
Any bialgebra comultiplication $\Delta :End(V)\to End(V)\otimes End(V)$ is defined by an isomorphism 
$F:V\otimes M\to V\otimes V$ (for some vector space $M$) such that 
\begin{equation}\label{pc}
F_{12}F_{13} \equiv F_{23}F_{12}\ mod(1\otimes Aut(M^{\otimes 2}))
\end{equation}
in the following way
$$\Delta_F (x) = F(x\otimes 1)F^{-1},\ \forall x\in End(V).$$
The bialgebras $(End(V), \Delta_F ),(End(V), \Delta_{F'})$ defined by the automorphisms 
$F,F'$ are isomorphic iff there exist an isomorphisms $f:V\to V, g:M\to M'$ such that 
\begin{equation}\label{re}
(f\otimes f)F = F'(f\otimes g).
\end{equation}
\end{th}
Proof:

Denoting $Hom_{End(V)}(V,V^{\otimes 2})$ by $M$ we get desired isomorphism (evaluation map):
$$F:V\otimes M\to V\otimes V$$
which by lemma \ref{hea} conjugates given inclusion $\Delta :End(V)\to End(V)\otimes End(V)$ with the standard one $End(V)\to End(V)\otimes End(M)$  
($x\mapsto (x\otimes 1)$) so that $\Delta (x) = \Delta_F (x) = F(x\otimes 1)F^{-1}$. 

Coassociativity of the comultiplication $\Delta_F$ is equivalent to the congruence (\ref{pc}).
Indeed, the equality between 
$$(\Delta_F\otimes I)\Delta_F (x) = F_{12}F_{13}(x\otimes 1\otimes 1)F_{13}^{-1}F_{12}^{-1}$$ 
and
$$(I\otimes \Delta_F )\Delta_F (x) = 
F_{23}F_{12}(x\otimes 1\otimes 1)F_{12}^{-1}F_{23}^{-1}$$
means that $F_{12}^{-1}F_{23}^{-1}F_{12}F_{13}$ lies in the centralizer of 
$End(V)\otimes 1\otimes 1$ in $End(V\otimes M\otimes M)$ which by lemma \ref{hea} coincides with 
$1\otimes End(M^{\otimes 2})$.
$\Box$

We say that the pair of isomorphisms $(F,\Phi), F:V\otimes M\to V\otimes V, \Phi\in Aut(M^{\otimes 2})$ satisfies to 
{\em modified pentagon equation} on vector spaces $V,M$ if
\begin{equation}\label{mpe}
F_{12}F_{13}\Phi_{23} = F_{23}F_{12}. 
\end{equation}
For any solution $(F,\Phi )$ of modified pentagon equation automorphism $F$ satisfies to the 
congruence (\ref{pc}) and vise versa 
any solution $F$ of the congruence (\ref{pc}) defines the solution $(F,\Phi )$ of modified 
pentagon equation by 
$$\Phi_{23} = F_{13}^{-1}F_{12}^{-1}F_{23}F_{12}.$$
Two solution $F,F'$ of the congruence \ref{pc} lie in the same equivalence class of the 
relation (\ref{re})
$$(f\otimes f)F = F'(f\otimes g),\ \mbox{ for some } f,g\in Aut(V)$$
if and only if for the pairs $(F,\Phi ),(F',\Phi' )$ we have 
$$(f\otimes f)F = F'(f\otimes g),\ (g\otimes g)\Phi = {\Phi}'(g\otimes g)\quad f\in Aut(V), g\in Aut(M).$$
Indeed,
$$(f\otimes f\otimes f)F_{23}F_{12} = F_{23}'(f\otimes f\otimes g)F_{12} = 
F_{23}'F_{12}'(f\otimes g\otimes g)$$
coincides with
$$(f\otimes f\otimes f)F_{12}F_{13}\Phi_{23} = F_{12}'(f\otimes g\otimes f)F_{13}\Phi_{23} = 
F_{12}'F_{13}'(f\otimes g\otimes g)\Phi_{23}$$
which means that $(g\otimes g)\Phi = {\Phi}'(g\otimes g)$. 
\begin{prop}\label{moppe}
If the pair $F,\Phi$ satisfies to modified pentagon equation (\ref{mpe}) 
then $\Phi$ is a solution to {\em pentagon equation}
\begin{equation}\label{pe}
\Phi_{12}\Phi_{13}\Phi_{23} = \Phi_{23}\Phi_{12}.
\end{equation}
\end{prop}
Proof:

Since $F_{14}\Phi_{24}F_{12}^{-1} = F_{12}^{-1}F_{24}$ the expression 
$$\Phi_{23}\Phi_{24}\Phi_{34} = 
F_{14}^{-1}\Phi_{23}F_{14}\Phi_{24}F_{12}^{-1}\Phi_{34}F_{12}$$
coincides with 
$$F_{14}^{-1}\Phi_{23}F_{12}^{-1}F_{24}\Phi_{34}F_{12} = 
F_{14}^{-1}\Phi_{23}F_{12}^{-1}F_{23}^{-1}F_{23}F_{24}\Phi_{34}F_{12}.$$
Which again using modified pentagon equations $\Phi_{23}F_{12}^{-1}F_{23}^{-1} =
F_{13}^{-1}F_{12}^{-1}$ and $F_{23}F_{24}\Phi_{34} = F_{34}F_{23}$ can be transformed into 
$$F_{14}^{-1}F_{13}^{-1}F_{12}^{-1}F_{34}F_{23}F_{12} = 
F_{14}^{-1}F_{13}^{-1}F_{34}F_{12}^{-1}F_{23}F_{12} = 
F_{14}^{-1}F_{13}^{-1}F_{34}F_{13}F_{13}^{-1}F_{12}^{-1}F_{23}F_{12},$$
which in his turn coincides with 
$$\Phi_{23}\Phi_{12}$$
since $F_{14}^{-1}F_{13}^{-1}F_{34}F_{13} = \Phi_{23}$ and 
$F_{13}^{-1}F_{12}^{-1}F_{23}F_{12} = \Phi_{12}$.
$\Box$

Note that the pair of identity operators $(I,I)$ is a solution to modified pentagon equation. 
\begin{exam}
Let $\Phi$ be a solution of pentagon equation. Then $(\Phi ,\Phi )$ is a solution to the modified pentagon equation. 
\end{exam}
The modified pentagon equation for the pair $(\Phi ,\Phi )$ is an ordinary pentagon equation for $\Phi$.$\Box$

Denote by $t = t_V:V\otimes V\to V\otimes V$ the permutation of tensor factors.
If $\Delta$ is a coproduct on the algebra $E$, then $t_E\Delta$ is also a coproduct. Applying this to endomorphism algebra $E=End(V)$ we get the following operation on solutions to modified pentagon equation.
\begin{exam}
Let $(F,\Phi )$ be a solution of modified pentagon equation. Then $(tF,\Phi^{-1}_{21})$ is also a solution.
In particular, for any solution $\Phi$ to pentagon equation $\Phi^{-1}_{21}$ is also a solution.
\end{exam}
Let us give the direct verification of this fact. 
First, conjugation modified pentagon equation for $(F,\Phi)$ with $t_{23}$ we get the following
$$F_{13}F_{12}\Phi_{32} = F_{32}F_{13}\quad \mbox{or}\quad F_{32}F_{13}\Phi_{32}^{-1} = F_{13}F_{12}.$$
Then, since $t_{12}t_{13} = t_{23}t_{12}$ the left side of modified pentagon equation for $(tF,\Phi^{-1}_{21})$
$$t_{12}F_{12}t_{13}F_{13}\Phi_{32}^{-1} = t_{12}t_{13}F_{32}F_{13}\Phi_{32}^{-1}$$
coincides with the right side
$$t_{23}F_{23}t_{12}F_{12} = t_{23}t_{12}F_{13}F_{12}.$$
$\Box$

Starting with the pair $(\Phi^{-1}_{21},\Phi^{-1}_{21})$ we get that $(\Phi^{-1}t,\Phi )$ is a solution to the modified pentagon equation for any solution of pentagon equation $\Phi$. In particular, the pair $(t,I)$ is a solution to modified pentagon equation. 

Tensor product of bialgebras $(End(V),\Delta_F),\ (End(V'),\Delta_{F'})$ corresponds to the following operation on solutions of modified pentagon equation. $\Box$
\begin{exam}
Let $(F,\Phi), (\Phi',F')$ be a solutions of pentagon equation on $(V,M)$ and $(V',M')$ respectively. 
Then $(t_{23}(F\otimes\F')t_{23},t_{23}(\Phi\otimes\Phi')t_{23})$ is a solution of pentagon equation on $(V\otimes V',M\otimes M')$.

In particular, for identity solution $(F,\Phi')=(I,I)$ the pair $(t_{23}(F\otimes I)t_{23},t_{23}(\Phi\otimes I)t_{23})$ 
is a solution of pentagon equation.
\end{exam}
As a partial case of previous example we have the following tensor product operation for solutions of pentagon equation.
\begin{exam}
Let $\Phi\in Aut(V^{\otimes 2}), \Phi'\in Aut({V'}^{\otimes 2})$ be a solutions of pentagon equation. 
Then $t_{23}(\Phi\otimes\Phi')t_{23}\in Aut((V\otimes V')^{\otimes 2})$ is a solution of pentagon equation.

In particular, for identity solution $\Phi'$ the operator $t_{23}(\Phi\otimes I)t_{23}$ 
is a solution of pentagon equation.
\end{exam}

\section{Hopf modules and examples of solutions to pentagon equation}\label{esmpe}
Here we give a way of constructing solutions to pentagon equation using Hopf modules.
In the slightly different form this construction appeared in the papers \cite{mil1,mil2}.

Let $H$ be a Hopf algebra. {\em Hopf $H$-module} \cite{mon} is a vector space $M$ with a structure of left $H$-module
$$\mu_M :H\otimes M\to M,\quad \mu_M (h\otimes m) = hm$$
and a structure of right $H$-comodule
$$\Delta_M :M\to M\otimes H,$$
compatible in the following sense:
$$\Delta_M (hm) = \Delta_H (h)\Delta_M (m)\quad (H\mbox{-linearity}).$$
\begin{prop}\label{homosol}
Let $H$ be a Hopf algebra and $M$ be a Hopf $H$-module with multiplication $\mu$ and comultiplication $\Delta$. 
Then the formula
$$\Phi_M = (I\otimes\mu )(\Delta\otimes I)$$
defines a solution to pentagon equation on vector space $M\otimes M$.
\end{prop}
Proof:

We will use Sweedler's notation for comodule structure $\Delta_M (m) = \sum_{(m)}m_{(0)}\otimes m_{(1)}$. 
For example, the map $\Phi_M$ has a form: 
$$\Phi(m\otimes n) = \sum_{(m)}m_{(0)}\otimes m_{(1)}n.$$
It can be checked directly that the inverse map is given by the formula:
$$\Phi^{-1}(m\otimes n) = \sum_{(m)}m_{(0)}\otimes S(m_{(1)})n.$$
For the sake of simplicity we will omit summation signs in verification of pentagon equation:
$$\Phi_{12}\Phi_{13}\Phi_{23}(m\otimes n\otimes l) = \Phi_{12}\Phi_{13}(m\otimes n_{(0)}\otimes n_{(1)}l) = $$
$$\Phi_{12}(m_{(0)}\otimes n_{(0)}\otimes m_{(1)}n_{(1)}l) = m_{(0)}\otimes m_{(1)}n_{(0)}\otimes m_{(2)}n_{(1)}l.$$
On the other side 
$$\Phi_{23}\Phi_{12}(m\otimes n\otimes l) = \Phi_{23}(m_{(0)}\otimes m_{(1)}n\otimes l) = 
m_{(0)}\otimes m_{(1)}n_{(0)}\otimes m_{(2)}n_{(1)}l.$$
$\Box$
\begin{exam}\label{ehm}
The Hopf algebra $H$ considered as a left module and a right comodule over itself is a Hopf module.
We shall call this Hopf module {\em trivial}. 
\end{exam}
We can define Hopf module structure on any multiple $M\otimes V$ of Hopf module $M$. 
It can be checked directly that corresponding solutions of pentagon equations are related as follows:
$$\Phi_{M\otimes V} = t_{23}(\Phi_M\otimes V)t_{23}.$$

The next theorem is known as "Fundamental theorem of Hopf modules" (see, for example, \cite{mon}). 
It states that up to a multiplication by vector space any Hopf module is isomorphic to 
trivial one. Here by a homomorphism of Hopf modules we mean a homomorphism preserving module and comodule structures. 
For a Hopf module $M$ denote by $M_H = \{ m\in M, \Delta_M (m)=m\otimes 1\}$ the subspace of {\em coinvariant} elements. 
\begin{th}
Let $M$  be a Hopf $H$ module. Then the map
$$H\otimes M_H\to M,\quad h\otimes m\mapsto hm$$
is an isomorphism of Hopf modules where Hopf module structure on $H\otimes M_H$ is induced by the natural
Hopf module structure on $H$ (see example \ref{ehm}). 
\end{th}
Hint of the proof:

Construct inverse map $M\to H\otimes M_H$ as $m\mapsto \sum_{(m)}m_{(2)}\otimes S^{-1}(m_{(1)})m_{(0)}$ where
$S^{-1}$ is inverse map to the antipode $S$ of the Hopf algebra $H$.
$\Box$

\begin{cor}\label{hms}
For any Hopf $H$-module $M$ the solution of pentagon equation $\Phi_M$ is isomorphic to $\Phi_H\otimes M_H$. 
\end{cor}

\section{Construction of solutions for modified pentagon equation}\label{esm}

Here we give a way of constructing solutions to modified pentagon equation for a given solution of pentagon equation $\Phi$  
corresponding to Hopf $H$-module $M$. The initial data is more peculiar then in the case of pentagon equation and consists of:

right $H$-{\em module coalgebra} $L$ with structure maps
$$\Delta_L :L\to L\otimes L,\quad (\mbox{coalgebra comultiplication}),$$
$$\mu_L:L\otimes H\to L,\quad (H-\mbox{module structure});$$
right $L$-comodule $V$
$$\Delta_V:V\to V\otimes L$$
and the pairing
\begin{equation}\label{pai}
\pi :L\otimes M\to V.
\end{equation}
In addition to ordinary axioms of module coalgebra:
$$(\Delta_L\otimes I)\Delta_L = (I\otimes\Delta_L)\Delta_L:L\to L\otimes L\otimes L,$$
$$\mu_L(\mu_L\otimes I) = \mu_L(I\otimes\mu_H):L\otimes H\otimes H\to L,$$
$$\Delta_L\mu_L = (\mu_L\otimes\mu_H)t_{23}(\Delta_L\otimes\Delta_H):L\otimes H\to L\otimes L$$
and of comodule
$$(\Delta_V\otimes I)\Delta_V = (I\otimes\Delta_L)\Delta_V:V\to V\otimes L\otimes L$$
we need following properties for the pairing \ref{pai}:
$$\pi (\mu_L\otimes I) = \pi (I\otimes\mu_M):L\otimes H\otimes M\to V,$$
$$\Delta_V\pi = (\pi\otimes\mu_L)t_{23}(\Delta_L\otimes\Delta_M):L\otimes M\to V\otimes L,$$
where $\mu_M:H\otimes M\to M$, $\Delta_M:M\to M\otimes H$ are $H$-module and comodule structure on $M$ respectively. 
Note that the first property means that $\pi$ factors through tensor product over $H$ $\pi:L\otimes M\to L\otimes_H M\to V$. Since $M\simeq H\otimes M^H$ as an $H$-module $\pi$is defined by the map $\overline\pi :L\otimes M^H\to V$. The second property says that $\overline\pi$ is a homomorphism of $L$-comodules where $L$-comodule structure on $L\otimes M^H$ comes from those on $L$. 

Having all of these we can define a map $F_V:V\otimes M\to V\otimes V$ by 
$F_V = (I\otimes\pi )(\Delta_V\otimes I)$. 
\begin{prop}
The map $F_V$ satisfies to modified pentagon equation $F_{12}F_{13}\Phi_{23} = F_{23}F_{12}$, where $\Phi$ is the solution of pentagon equation corresponding to Hopf $H$-module $M$. 
\end{prop}
Proof:

The proof is similar to the proof of proposition \ref{homosol}. $\Box$

Now we give some additional data providing invertibility of the map $F_V$. 
\begin{prop}
Suppose that our coalgebra is equipped with the counite $\varepsilon_L:L\to k$ and we are given by a pairing $\nu:L\otimes V\to M$ such that the compositions
$$\nu (I\otimes\pi )(\Delta_L\otimes I):L\otimes M\to M,$$
$$\pi (I\otimes\nu )(\Delta_L\otimes I):L\otimes V\to V$$
have a form $\varepsilon_L\otimes I$. Then the map $F_V$ is invertible with the inverse $F_V^{-1} = (I\otimes\nu )(\Delta_L\otimes I)$.
\end{prop}
Proof:

Consists of direct checking.$\Box$ 

Consider as an example the case when vector space $M^H$ is one dimensional. 
\newline
Let $L$ be a right $H$-module coalgebra . Define
$$F_L = (I\otimes\mu_L )(\Delta_L\otimes I):L\otimes H\to L\otimes L,$$
where $\mu_L :L\otimes H\to L$ is the action map and $\Delta_L :L\to L\otimes L$ is the 
coalgebra structure. 
Right $H$-module coalgebra $L$ is called {\em Galois} \cite{mon} if the map $F_L$ is an isomorphism. 
\begin{exam}\label{galmpe}
Then the pair $(F_L ,\Phi_H )$ is a solution of modified pentagon equation. 
\end{exam}

\section{Solutions of pentagon equation}\label{sspe}
The theorem of Baaj and Skandalis (Theorem 4.10. of \cite{sask}) says that any finite 
dimensional unitary solution of pentagon equation corresponds to some Hopf module over 
some Hopf $C^*$-algebra. In this section we give a proof this theorem without any unitary 
assumptions which works over arbitrary field. We basically follow the scheme of \cite{sask} 
with some minor changes. Namely we use Hamilton-Cayley's theorem defining counite and 
antipode and the notion of Hopf module in formulation of the theorem. 

In papers of G.Militaru \cite{mil1,mil2} the variant of the theorem of Baaj and Skandalis was 
proved for more general situation of arbitrary (not necessary invertible) solution of 
pentagon equation but the resulting Hopf module there is a Hopf module over bialgebra which 
does not posses an antipode in general.

According to the theorem \ref{cea} two solutions $(\Phi ,\Phi )$ and $(\Phi^{-1}t,\Phi )$ of 
the modified pentagon equation define two comultiplications on the endomorphism algebra 
$End(M)$
$$\Delta_l (x) = \Delta_{\Phi}(x) = \Phi (x\otimes 1)\Phi^{-1},$$
$$\Delta_r (x) = \Delta_{\Phi^{-1}t}(x) = \Phi^{-1}(1\otimes x)\Phi ,$$
which give two multiplications $\mu_l ,\mu_r$ on the dual space $End(M)^*$. 
\begin{lem}
The maps defined by 
$$\lambda (\omega ) = (\omega\otimes I)(\Phi ),\ \rho (\omega ) = (I\otimes\omega )(\Phi )$$
are a homomorphisms of bialgebras
$$\lambda :(End(M)^* ,\mu_r ,\Delta )\to (End(M),\cdot ,\Delta_l ),\ \rho :(End(M)^* ,\mu_l ,\Delta )\to (End(M),\cdot ,\Delta_r ),$$
where $\cdot$ is the composition of endomorphisms in $End(M)$ and $\Delta$ denotes a comultiplication on $End(M)^*$ dual to $\cdot$. 
\end{lem}
Proof:

The homomorphism property for $\lambda$ follows from pentagon equation for $\Phi$:
$$\lambda (\omega )\lambda (\omega') = (\omega\otimes\omega'\otimes I)(\Phi_{13}\Phi_{23}) =
(\omega\otimes\omega'\otimes I)(\Phi_{12}^{-1}\Phi_{23}\Phi_{12}) = 
\lambda (\omega *_r\omega').$$
Analogously for $\rho$:
$$\rho (\omega )\rho (\omega') = (I\otimes\omega\otimes\omega')(\Phi_{12}\Phi_{13}) =
(I\otimes\omega\otimes\omega')(\Phi_{23}\Phi_{12}\Phi_{23}^{-1}) = 
\rho (\omega *_l\omega').$$
Now check that $\lambda$ is a homomorphism of coalgebras, that is $\Delta_l\lambda = (\lambda\otimes\lambda )\Delta$:
$$(\lambda\otimes\lambda )(\Delta (\omega )) = (\Delta (\omega)\otimes I\otimes I)(\Phi_{13}\Phi_{24}) = (\omega\otimes I\otimes I)(\Phi_{12}\Phi_{13}) = $$
$$(\omega\otimes I\otimes I)(\Phi_{23}\Phi_{12}\Phi_{23}^{-1}) = \Delta_l\lambda (\omega ).$$
Analogously, for $\Delta_r\rho = (\rho\otimes\rho )\Delta_r$:
$$(\rho\otimes\rho )(\Delta (\omega )) = (I\otimes I\otimes\Delta_r (\omega ))(\Phi_{13}\Phi_{24}) = (I\otimes I\otimes\omega )(\Phi_{13}\Phi_{23}) = $$
$$(I\otimes I\otimes\omega )(\Phi_{12}^{-1}\Phi_{23}\Phi_{12}) = \Delta_r\rho (\omega ).$$
$\Box$

In particular, $im(\lambda ),\ im(\rho )$ are a bialgebras and subalgebras of $End(M)$. 
By the definition $\Phi\in im(\rho )\otimes im(\lambda )\subset End(M)$ and formula
$$(\rho (\omega ),\lambda (\omega')) = (\omega'\otimes\omega )(\Phi )$$
defines non-degenerated pairing $(\ ,\ ):im(\rho )\otimes im(\lambda )\to k$. 
\begin{lem}
The pairing $(\ ,\ )$ is a pairing of bialgebras, e.g.
$$(\Delta_r (\rho (\omega )),\lambda (\omega')\otimes\lambda (\omega'')) = (\rho (\omega ),\lambda (\omega')\lambda (\omega''))$$
and 
$$(\rho (\omega )\otimes\rho (\omega'),\Delta_l (\lambda (\omega'')) = (\rho (\omega )\rho (\omega'),\lambda (\omega'')).$$ 
\end{lem}
Proof:

The verification is direct:
$$(\Delta_r (\rho (\omega )),\lambda (\omega')\otimes\lambda (\omega'')) = ((\rho\otimes\rho )(\Delta (\omega )),\lambda (\omega')\otimes\lambda (\omega'')) = $$
$$(\Delta (\omega )\otimes\omega'\otimes\omega'')(\Phi_{13}\Phi_{24}) = (\omega\otimes\omega'\otimes\omega'')(\Phi_{12}\Phi_{13}) = (\omega\otimes\omega'\otimes\omega'')(\Phi_{23}\Phi_{12}\Phi_{23}^{-1}) = $$
$$(\omega\otimes\omega' *_l\omega'')(\Phi ) = (\rho (\omega ),\lambda (\omega' *_l\omega'')) = (\rho (\omega ),\lambda (\omega')\lambda (\omega''))$$
and
$$(\rho (\omega )\otimes\rho (\omega'),\Delta_l (\lambda (\omega''))) = (\rho (\omega )\otimes\rho (\omega'),(\lambda\otimes\lambda )(\Delta (\omega''))) = $$ 
$$(\omega\otimes\omega'\otimes\Delta (\omega''))(\Phi_{13}\Phi_{24}) = (\omega\otimes\omega'\otimes\omega'')(\Phi_{13}\Phi_{23}) = (\omega\otimes\omega'\otimes\omega'')(\Phi_{12}^{-1}\Phi_{12}\Phi_{12}) = $$
$$(\omega *_r\omega'\otimes\omega'')(\Phi ) = (\rho (\omega *_r\omega'),\lambda (\omega'')) = (\rho (\omega )\rho (\omega'),\lambda (\omega'')).$$
$\Box$

\begin{lem}
Subalgebras $im(\lambda ),\ im(\rho )\subset End(M)$ are unital, e.g. 
endomorphism $I$ belongs to $im(\rho)$ and $im(\lambda )$.
\end{lem}
Proof:

Here we use the fact that $M$ is finite dimensional vector space. 

Let $\chi_\Phi (t)\in k[t]$ be the characteristic polynomial of $\Phi\in End(M^{\otimes 2})$. 
Since $\Phi$ is invertible $\chi_\Phi (0)\not= 0$. By Hamilton-Cayley's theorem $\chi_\Phi (\Phi )=0$. Hence we can represent identity as a linear combinations of powers of $\Phi$. Namely, 
$I = f(\Phi )\in im(\rho)\otimes im(\lambda )$ where $f(t) = det(\Phi )^{-1}\chi_\Phi (t) - 1$. 
$\Box$

Now we construct a counite for the bialgebra $im(\lambda )$ using duality between $im(\lambda )$ and $im(\rho )$. 
\begin{prop}
Let $\varepsilon\in\rho^{-1}(I)$, that is $\rho (\varepsilon ) = I$. Then the restriction of $\varepsilon$ to $im(\lambda )$ is a homomorphism of algebras $\varepsilon :im(\lambda )\to k$. 

This homomorphism is a counite with respect to the coproduct $\Delta_l$, that is
$$(\varepsilon\otimes I)\Delta_l = (I\otimes\varepsilon )\Delta_l = I.$$
\end{prop}
Proof:

Note that for any $\omega\in End(M)^*$ 
$$\varepsilon (\lambda (\omega )) = (\omega\otimes\varepsilon )(\Phi ) = \omega (I)$$
and  for any $\omega ,\omega'\in End(M)^*$ 
$$(\omega *_r\omega' )(I) = (\omega\otimes\omega' )(\Phi^{-1}(I\otimes I)\Phi ) = 
(\omega\otimes\omega' )(I\otimes I) = \omega (I)\omega' (I).$$ 
Now 
$$\varepsilon (\lambda (\omega )\lambda (\omega' )) = \varepsilon (\lambda (\omega *_r\omega' )) = 
(\omega *_r\omega' )(I) = \omega (I)\omega' (I) = \varepsilon (\lambda (\omega ))\varepsilon (\lambda (\omega' )).$$
The counite axiom can be verified directly:
$$(\varepsilon\otimes I)\Delta (\lambda (\omega )) = (\varepsilon\otimes I)(\lambda\otimes\lambda )(\Delta (\omega )) = 
(\omega\otimes\varepsilon\otimes I)(\Phi_{12}\Phi_{13}) = (\omega\otimes I)(\Phi ) = \lambda (\omega ).$$
Analogously,
$$(I\otimes\varepsilon )\Delta (\lambda (\omega )) = (I\otimes\varepsilon )(\lambda\otimes\lambda )(\Delta (\omega )) = 
(\omega\otimes I\otimes\varepsilon )(\Phi_{12}\Phi_{13}) = (\omega\otimes I)(\Phi ) = \lambda (\omega ).$$
$\Box$

It follows by the uniqueness of a counite that all $\varepsilon\in\rho^{-1}(I)$ give the same restriction to $im(\lambda)$. 

The next step is to define the antipode of the bialgebra $im(\lambda )$.
\begin{prop}
The map $S:im(\lambda )\to im(\lambda )$ given by $S(\lambda (\omega)) = \sigma(\omega) = (\omega\otimes I)(\Phi^{-1})$ is an antipode for the bialgebra $im(\lambda )$.
\end{prop}
Proof:

First we need to check the correctness of the definition. Note that $\omega(im(\rho ))=0$ for $\omega\in ker(\lambda )$. 
Indeed,
$$\omega(\rho (\omega')) = (\omega\otimes\omega')(\Phi) = 0.$$
Since $\Phi^{-1}$ belongs to $im(\rho )\otimes im(\lambda )$ we have that $(\omega\otimes I)(\Phi^{-1}) = 0$ for any $\omega\in ker(\lambda )$. This proves that the assertion $S(\lambda (\omega)) = (\omega\otimes I)(\Phi^{-1})$ defines a map from $im(\lambda )$. To check that this is a map to $im(\lambda )$ is enough to note that since $\Phi^{-1}\in im(\rho )\otimes im(\lambda )$ any expression of the form $(\omega\otimes I)(\Phi^{-1})$ lies in $im(\lambda )$. 

Now we are ready to verify antipode axioms:
$$(I\otimes S)\Delta_l(\lambda (\omega)) = (I\otimes S)(\lambda\otimes\lambda )(\Delta(\omega)) = (\lambda\otimes\sigma)(\Delta(\omega)) = $$
$$(\Delta(\omega)\otimes I\otimes I)(\Phi_{13}\Phi_{24}^{-1}) = (\omega\otimes I\otimes I)(\Phi_{12}\Phi_{13}^{-1})$$
so
$$\mu(I\otimes S)\Delta_l (\lambda (\omega)) =  (\omega\otimes I)(\Phi_{12}\Phi_{12}^{-1}) = \omega(I)I.$$
Analogously, 
$$(S\otimes I)\Delta_l(\lambda (\omega)) = (S\otimes I)(\lambda\otimes\lambda )(\Delta(\omega)) = (\sigma\otimes\lambda )(\Delta(\omega)) = $$
$$(\Delta(\omega)\otimes I\otimes I)(\Phi_{13}^{-1}\Phi_{24}) = (\omega\otimes I\otimes I)(\Phi_{12}^{-1}\Phi_{13})$$
so
$$\mu(S\otimes I)\Delta_l (\lambda (\omega)) =  (\omega\otimes I)(\Phi_{12}^{-1}\Phi_{12}) = \omega(I)I.$$
$\Box$

We have proven that $H = im(\lambda)$ is a Hopf algebra. 
By the construction $M$ is a left module over $H$ and over dual Hopf algebra $H^*\simeq im(\rho)$. 
In the next proposition we construct $H$-comodule structure on $M$ which corresponds to $H^*$-module structure
and prove that this structure is $H$-linear, so that $M$ is a Hopf $H$-module. 
\begin{prop}
The formula $r_\Phi (m) = \Phi (m\otimes I)$ defines right $H$-comodule structure on $M$.

This comodule structure satisfies to the equation:
$$r_\Phi (\lambda(\omega)m) = \Delta_l (\lambda(\omega))r_\Phi (m).$$
\end{prop}
Proof:

Since $\Phi$ lies in $im(\rho)\otimes im(\lambda)$ then $\Phi (v\otimes I)$ belongs to $M\otimes im(\lambda)$. 
Coassociativity of comodule structure is a consequence of pentagon equation:
$$(r_\Phi\otimes I)r_\Phi (m) = (\Phi_{12}\Phi_{13})(m\otimes I\otimes I) = (\Phi_{23}\Phi_{12}\Phi_{23}^{-1})(m\otimes I\otimes I) = (I\otimes\Delta_H )r_\Phi (m).$$
Transforming left side of $H$-linearity equation:
$$r_\Phi (\lambda(\omega)m) = (\omega\otimes I\otimes I)(\Phi_{23}\Phi_{12})(I\otimes m\otimes I) = 
(\omega\otimes I\otimes I)(\Phi_{12}\Phi_{13}\Phi_{23})(I\otimes m\otimes I)$$
we can see that it coincides with $\Delta_l (\lambda(\omega))r_\Phi (m)$. 
$\Box$

Now we are ready to prove the main theorem of this section (see theorem 4.10 from \cite{sask}) which says that up to multiplication by a vector space any solution of pentagon equation corresponds to some Hopf algebra. 
\begin{th}\label{spe}
The solution to the pentagon equation $\Phi_M$ corresponding to the Hopf $H$-module structure on $M$ coincides with $\Phi$. 

In particular, the multiplication map $H\otimes M_H\to M$ is an isomorphism and identifies $\Phi$ with $\Phi_H\otimes M_H$,
where $M_H = \{ m\in M,\ \Phi(m\otimes n) = m\otimes n\ \forall n\in M\}$. 
\end{th}
Proof:

By the definition of Hopf $H$-module structure on $M$: 
$$\Phi_M (m\otimes n) = (I\otimes\mu_M )(\Delta_M\otimes I)(m\otimes n) = \Phi(m\otimes n).$$
Now by the corollary \ref{hms} the multiplication map $H\otimes M_H\to M$ is an isomorphism and identifies $\Phi$ with $\Phi_H\otimes M_H$, where $M_H$ is a space of coinvariants of the comodule structure $\Delta_M$. Using the definition of $\Delta_M$ we can see that $M_H$ coincides with the space $\{ m\in M,\ \Phi(m\otimes n) = m\otimes n\ \forall n\in M\}$.
$\Box$

\section{Modified pentagon equation}

In this section we apply methods of the previous section to solutions of modified pentagon equation.

\begin{prop}
Let $(F,\Phi )$ be a solution of modified pentagon equation on $(V,M)$. Then the map 
$\Delta_{F,\Phi}:Hom(M,V)\to Hom(M,V)\otimes Hom(M,V)$ defined by
$$\Delta_{F,\Phi}(x) = F(x\otimes 1)\Phi^{-1}$$
is a coassociative comultiplication. 

The map 
$$\lambda_F:(End(V)^*,\Delta)\to (Hom(M,V),\Delta_{F,\Phi})$$
is a homomorphism of coalgebras and the map 
$$\rho_F:(Hom(M,V)^* ,\mu_{F,\Phi})\to End(V),\ \rho_F(\omega) = (I\otimes\omega)(F)$$
is a homomorphism of algebras, where $\mu_{F,\Phi}$ is a multiplication on the dual space $Hom(M,V)^*$ 
defined by $\Delta_{F,\Phi}$. 
\end{prop}
Proof:

The coassociativity of $\Delta_{F,\Phi}$ is a direct consequence of modified pentagon equation for 
$(F,\Phi )$ and pentagon equation for $\Phi$ which being combined together state that 
$$F_{12}^{-1}F_{23}^{-1}F_{12}F_{13} = \Phi_{23}^{-1} = 
\Phi_{12}^{-1}\Phi_{23}^{-1}\Phi_{12}\Phi_{13}.$$
This implies the coincidence between 
$$(\Delta_{F,\Phi}\otimes I)\Delta_{F,\Phi}(x) = 
F_{12}F_{13}(x\otimes 1\otimes 1)\Phi_{13}^{-1}\Phi_{12}^{-1}$$ 
and
$$(I\otimes \Delta_{F,\Phi})\Delta_{F,\Phi}(x) = 
F_{23}F_{12}(x\otimes 1\otimes 1)\Phi_{12}^{-1}\Phi_{23}^{-1}.$$
The homomorphism properties for $\lambda_F, \rho_F$ also follow from the modified pentagon equation
$$(\lambda_F\otimes\lambda_F)(\Delta (\omega )) = (\Delta (\omega)\otimes I\otimes I)(F_{13}F_{24}) = (\omega\otimes I\otimes I)(F_{12}F_{13}) = $$
$$(\omega\otimes I\otimes I)(F_{23}F_{12}\Phi_{23}^{-1}) = \Delta_{F,\Phi}\lambda (\omega ).$$
and
$$\rho_F(\omega )\rho_F(\omega') = (I\otimes\omega\otimes\omega')(F_{12}F_{13}) =$$
$$(I\otimes\omega\otimes\omega')(F_{23}F_{12}\Phi_{23}^{-1}) = \rho_F(\omega *_{F,\Phi}\omega').$$
$\Box$

In particular, $im(\lambda_F)$ is a coalgebra and $im(\rho_F)$ is an algebra. In the next proposition we prove that $im(\lambda_F)$ is an $im(\lambda_\Phi )$-module coalgebra and $im(\rho_F)$ is an $im(\lambda_\Phi )$-comodule algebra. 

\begin{prop}
The subspace $im(\lambda_F)\subset Hom(M,V)$ is closed under right multiplication by elements from $im(\lambda_\Phi )\subset End(M)$. Moreover, this $im(\lambda_\Phi )$-module structure on $im(\lambda_F)$ is compatible with comultiplication:
$$\Delta_{F,\Phi}(\lambda_F (\omega)\lambda_\Phi (\omega')) = \Delta_{F,\Phi}(\lambda_F (\omega))\Delta_\Phi (\lambda_\Phi (\omega')).$$
Dually, comultiplication $\Delta_F$ sends $im(\rho_F)$ into $im(\rho_F)\otimes im(\rho_\Phi )$
which gives a right $im(\lambda_\Phi )$-comodule structure on the algebra $im(\rho_F)$.
\end{prop}
Proof:

$$\lambda_F (\omega)\lambda_\Phi (\omega') = (\omega\otimes\omega'\otimes I)(F_{13}\Phi_{23}) =
(\omega\otimes\omega'\otimes I)(F_{12}^{-1}F_{23}F_{12}) = 
\lambda (\omega *_{\overline{F}}\omega').$$
Compatibility with comultiplication follows from the equality 
$$\Delta_{F,\Phi}(xy) = \Delta_{F,\Phi}(x)\Delta_{\Phi}(y),\quad \forall x\in Hom(M,V), y\in End(M).$$ 
It follows from the equality 
$$(\rho_F\otimes\rho_\Phi )(\Delta (\omega )) = (I\otimes I\otimes\Delta (\omega ))(F_{13}\Phi_{24}) = (I\otimes I\otimes\omega )(F_{13}\Phi_{23}) = $$
$$(I\otimes I\otimes\omega )(F_{12}^{-1}F_{23}F_{12}) = \Delta_F\rho_F (\omega )$$
that $\Delta_F(im(\rho_F))\subset im(\rho_F)\otimes im(\rho_\Phi )$. Hence algebra homomorphism $\Delta_F$ defines $im(\lambda_\Phi )$-comodule structure. 
$\Box$

In the next proposition we construct $im(\lambda_F)$-comodule structure on $V$ which is dual to $im(\rho_F)$-module structure. 
\begin{prop}
The formula $r_F (v) = F(v\otimes I)\in V\otimes im(\lambda_F)$ defines right $im(\lambda_F)$-comodule structure on $V$.
\end{prop}
Proof:

Since $F$ lies in $im(\rho_F)\otimes im(\lambda_F)$ then $F(v\otimes I)$ belongs to $V\otimes im(\lambda_F)$. 
Coassociativity of comodule structure is a consequence of pentagon equation:
$$(r_F\otimes I)r_F (v) = (F_{12}F_{13})(v\otimes I\otimes I) = (F_{23}F_{12}\Phi_{23}^{-1})(v\otimes I\otimes I) = (I\otimes\Delta_{F,\Phi})r_F (v).$$
$\Box$

Since $im(\lambda_F)$ is a subspace of $Hom(M,V)$ we have a pairing 
$$\mu_F :im(\lambda_F)\otimes M\to V,\quad \lambda_F(\omega)\otimes m\mapsto\lambda_F(\omega)m.$$
By the definition it has the following associativity property:
$$\lambda_F(\omega)(\lambda_\Phi(\omega)m) = (\lambda_F(\omega)\lambda_\Phi(\omega))m,\ \forall \omega\in End(V)^*, m\in M.$$
Next proposition shows compatibility of this pairing
with the $im(\lambda_\Phi )$-comodule structure on $M$ and the $im(\lambda_F)$-comodule structure on $V$. 
\begin{prop}
The $im(\lambda_\Phi )$-comodule structure on $M$ and the $im(\lambda_F)$-comodule structure on $V$ are compatible in the following way
$$r_F(\lambda_F(\omega)m) = \Delta_{F,\Phi}(\lambda_F(\omega))r_\Phi (m),\quad \forall \omega\in End(V)^*, m\in M.$$
\end{prop}
Proof:

Indeed, 
$$r_F(\lambda_F(\omega)m) = (\omega\otimes I\otimes I)(F_{23}F_{12})(I\otimes m\otimes I),$$
which by modified pentagon equation coincides with 
$$(\omega\otimes I\otimes I)(F_{12}F_{13}\Phi_{23})(I\otimes m\otimes I) = (\Delta(\omega)\otimes I\otimes I)(F_{13}F_{24}\Phi_{34})(I\otimes m\otimes I) =$$
$$(\lambda_F\otimes\lambda_F)(\Delta(\omega))r_\Phi (m) = \Delta_{F,\Phi}(\lambda_F(\omega))r_\Phi (m).$$
$\Box$ 

Now we have all necessary properties for the pair $(V,M)$ and we are ready to prove main theorem of this section. 
\begin{th}
The solution to the modified pentagon equation $F_{(V,M)}$ corresponding to the pair $(V,M)$ (see section \ref{esm}) coincides with $F$.  
\end{th}
Proof:

By the definition of $F_{(V,M)}$: 
$$F_{(V,M)}(v\otimes m) = (I\otimes\mu_F)(r_F\otimes I)(v\otimes m) = F(v\otimes m).\Box$$

\section{Categorical reformulation}

We start with description of $k$-linear monoidal structures on the category of (finite dimensional) vector spaces $\V{\it ect}_k$ over $k$. The standard reference for the theory of monoidal categories and Tannaka-Krein theory is \cite{dm}. For the sake of completeness we give some basic definitions of the theory. 
A {\em monoidal category} is a category $\G$  with a bifunctor
$$\otimes :\G\times\G\to\G\qquad (X,Y) \mapsto X\otimes Y$$
which is called {\em tensor (or monoidal) product}. This functor is supposed to be equipped with a functorial collection of isomorphisms (so-called {\em associativity constraint})
$$\varphi_{X,Y,Z} : X\otimes(Y\otimes Z)\rightarrow (X\otimes Y)\otimes Z\qquad\forall X,Y,Z \in\G$$
satisfying the so-called {\em pentagon axiom}:

the diagram
\begin{equation}\label{pax}\begin{array}{ccccc}
X\otimes(Y\otimes(Z\otimes W)) &\stackrel{\varphi_{X,Y,Z\otimes W}}{\longrightarrow} & (X\otimes Y)\otimes(Z\otimes W) & \stackrel{\varphi_{X\otimes Y,Z,W}}{\longrightarrow} &  ((X\otimes Y)\otimes Z)\otimes W \\
\downarrow{\footnotesize I_X\otimes\varphi_{Y,Z,W}} & & & & {\footnotesize \varphi_{X,Y,Z}\otimes I_W}\uparrow \\
X\otimes((Y\otimes Z)\otimes W) & & \stackrel{\varphi_{X,Y\otimes Z,W}}{\longrightarrow} & & (X\otimes(Y\otimes Z))\otimes W 
\end{array}\end{equation}

is commutative for any objects $X,Y,Z,W \in {\cal G}$.

By the functoriality any tensor product functor which is just a $k$-bi-linear functor $\boxtimes:\V{\it ect}_k\times\V{\it ect}_k\to\V{\it ect}_k$ is defined by the vector space $k\boxtimes k = M$ ({\em multiplicity space}). Denote by $\boxtimes_M$ the tensor product functor corresponding to a given vector space $M$. We have $U_1\boxtimes_M U_2 = U_1\otimes M\otimes U_2$ for arbitrary vector spaces $U_i$. 

Now we establish the relation between associativity constraint for the tensor product $\boxtimes_M$ and pentagon equation. 
\begin{prop}
There is one-to-one correspondence between structures of monoidal category on the category of finite dimensional vector spaces and solution of pentagon equation.
\end{prop}
Proof:

By functoriality any associativity constraint $\varphi$ for $\boxtimes_M$ is given by the automorphism $\Phi\in Aut(M^{\otimes 2})$:
$$M\otimes M = k\boxtimes_M M = k\boxtimes_M (k\boxtimes_M k)\stackrel{\varphi_{k,k,k}}{\longrightarrow} (k\boxtimes_M k)\boxtimes_M k = M\boxtimes_M k = M\otimes M$$ 
and has the following general form   
$$\begin{array}{ccc}
U_1\boxtimes_M (U_2\boxtimes_M U_3) & 
\stackrel{\varphi_{U_1 , U_2 , U_3}}{\longrightarrow} & 
(U_1\boxtimes_M U_2)\boxtimes_M U_3 \\
||                    &                                 &                    || \\
U_1\boxtimes_M (U_2\otimes M\otimes U_3 ) &  & 
(U_1\otimes M\otimes U_2 )\boxtimes_M U_3  \\ 
||                    &                                 &                    || \\
U_1\otimes M\otimes U_2\otimes M\otimes U_3  & \stackrel{\Phi_{24}}{\longrightarrow} & 
U_1\otimes M\otimes U_2\otimes M\otimes U_3  \\ 
\end{array}$$ 
In particular, 
$$\varphi_{k\boxtimes_M k,k,k} = \varphi_{M,k,k} = \Phi_{23},$$
$$\varphi_{k,k\boxtimes_M k,k} = \varphi_{k,M,k} = \Phi_{13},$$
$$\varphi_{k,k,k\boxtimes_M k} = \varphi_{k,k,M} = \Phi_{12}.$$
Combining it with 
$$I_k\boxtimes_M\varphi_{k,k,k} = \Phi_{23},\quad\varphi_{k,k,k}\boxtimes_M I_k = \Phi_{12}$$
we see that pentagon axiom (\ref{pax}) for $\varphi$ is equivalent to the pentagon equation for $\Phi$
$$\Phi_{12}\Phi_{13}\Phi_{23} = \Phi_{23}\Phi_{12}.\Box$$

For a solution of pentagon equation $\Phi\in Aut(M^{\otimes 2})$ denote by $\G_\Phi$ the corresponding monoidal category structure on the category of finite dimensional vector spaces. 

Now we discuss the relation between $k$-linear monoidal functors from the category $\G_\Phi$ to vector spaces $\V{\it ect}_k$ and modified pentagon equation. 
A {\em monoidal functor} between monoidal categories $\G$
and $\cH$ is a functor $\omega:\G\to\cH$ , which is equipped with a functorial collection
of isomorphisms (the so-called {\em monoidal structure})
$$\omega_{X,Y}:\omega(X\otimes Y)\to\omega(X)\otimes\omega(Y)\quad\forall X,Y\in\G,$$
for which the following diagram is commutative for any objects $X,Y,Z \in {\cal G}$
\begin{equation}\label{dmf}\begin{array}{ccccc}
F(X\otimes(Y\otimes Z)) & \stackrel{\omega_{X,Y\otimes Z}}{\longrightarrow} & \omega(X)\otimes\omega(Y\otimes Z) & \stackrel{I\otimes\omega_{Y,Z}}{\longrightarrow}& \omega(X)\otimes(\omega(Y)\otimes\omega(Z)) \\
\downarrow{\footnotesize \omega(\varphi_{X,Y,Z})} & & & & \downarrow{\footnotesize \varphi_{\omega(X),\omega(Y),\omega(Z)}} \\
\omega((X\otimes Y)\otimes Z) & \stackrel{\omega_{X\otimes Y,Z}}{\longrightarrow} &
\omega(X\otimes Y)\otimes\omega(Z) & \stackrel{\omega_{X,Y}\otimes I}{\longrightarrow} &
(\omega(X)\otimes\omega(Y))\otimes\omega(Z) \\
\end{array}\end{equation}

\begin{prop}\label{mofu}
There is one-to-one correspondence between $k$-linear monoidal functors $\G_\Phi\to\V{\it ect}_k$ and vector spaces $V$ with isomorphisms $F:V\otimes M\to V\otimes V$ satisfying modified pentagon equation $F_{12}F_{13}\Phi_{23} = F_{23}F_{12}.$
\end{prop}
Proof:

By functoriality any $k$-linear functor $\omega$ is defined by its evaluation $\omega(k) = V$. For arbitrary vector space $U$ it gives $\omega(U) = V\otimes U$. 
Monoidal structure of such a functor is given by the isomorphism $F:V\otimes M\to V\otimes V$:
$$V\otimes M = \omega(M) = \omega(k\boxtimes_M k)\stackrel{\omega_{k,k}}{\longrightarrow}\omega(k)\otimes\omega(k) = V\otimes V$$
and has the following general form 
$$\begin{array}{ccc}
\omega(U_1\boxtimes_M U_2 ) & 
\stackrel{\omega_{U_1 ,U_2}}{\longrightarrow} & 
\omega(U_1 )\otimes\omega(U_2 ) \\
||                    &                                 &                    || \\
\omega(U_1\otimes M\otimes U_2 ) &  & 
\omega(U_1 )\otimes \omega(U_2 ) \\ 
||                    &                                 &                    || \\
W\otimes U_1\otimes M\otimes U_2 & \stackrel{F_{13}}{\longrightarrow} & 
W\otimes U_1\otimes W\otimes U_2 \\ 
\end{array}$$ 
In particular, 
$$\omega_{k\boxtimes_M k,k} = \omega_{M,k} = F_{13},\quad \omega_{k,k\boxtimes_M k} = \omega_{k,M} = F_{12}.$$
Combining it with 
$$I_{\omega(k)}\otimes \omega_{k,k} = I_V\otimes F = \omega_{23},$$
$$\omega_{k,k}\otimes I_{\omega(k)} = F\otimes I_V = F_{12}$$
and with 
$$\omega(\phi_{k,k,k}) = \omega(\Phi ) = \Phi_{23}$$
we can see that compatibility axiom for monoidal structure (\ref{dmf}) is equivalent to the modified 
pentagon equation for $(F,\Phi )$
$$F_{12}F_{13}\Phi_{23} = F_{23}F_{12}.\Box$$

Denote by $\omega_F$ the functor $\G_\Phi\to\V{\it ect}_k$. General Tannaka-Krein theory \cite{dm} predicts existence of bialgebra structure on the algebra of natural transformations $End(\omega)$ of any monoidal functor $\omega$ from a semisimple category to the category of vector spaces. It is not hard to verify that $End(\omega_F)$ coincides with the bialgebra $(End(V),\Delta_F)$ defined in the first section. 

In general, we say that the pair $(V,F)$ consisting of an object $V$ of some $k$-linear monoidal category $\G$ and an isomorphism $v:V\otimes M\to V\otimes V$ (for some vector space $M$) is an {\em idempotent object} in $\G$ if $v$ satisfies to modified pentagon equation 
$$\begin{array}{ccccc}
V\otimes M\otimes M & \stackrel{v\otimes I}{\longrightarrow} & V\otimes V\otimes M & \stackrel{I\otimes v}{\longrightarrow} & V\otimes V\otimes V \\
\downarrow {\footnotesize I\otimes\Phi} & & & & \\
V\otimes M\otimes M & & & & {\footnotesize v\otimes I}\uparrow \\
{\footnotesize \downarrow c_{V,M}}\otimes I & & & & \\
M\otimes V\otimes M & \stackrel{I\otimes v}{\longrightarrow} & M\otimes V\otimes V & \stackrel{c_{M,V}\otimes I}{\longrightarrow} & V\otimes M\otimes V \\
\end{array}$$
or $v_{12}v_{13}\Phi_{23} = v_{23}v_{12}$ for some automorphism $\Phi\in Aut(M^{\otimes 2})$. For simplicity we omit brackets and associativity isomorphism in category $\G$.  
Then automatically (see proposition \ref{moppe}) $\Phi$ is a solution to pentagon equation. 

The examples of such objects are Hopf modules (see section \ref{esmpe} for definition). To clarify this we need the categorical characterization of Hopf modules. 

Let $\G$ be a $k$-linear monoidal category and $\omega:\G\to\V{\it ect}_k$ be a $k$-linear monoidal functor.  
We call the object $E\in\C$ {\em integral} with respect to $\omega$ if it is equipped with functorial collection of 
isomorphisms $e_{X}:E\otimes\omega(X)\to E\otimes X$ such that $e_1 = I$ and 
\begin{equation}\label{dinto}\begin{array}{ccc}
E\otimes\omega(X\otimes Y)        & \stackrel{e_{X\otimes Y}}{\longrightarrow} & E\otimes X\otimes Y \\
\downarrow {\footnotesize I\otimes\omega_{X,Y}} &                 & {\footnotesize e_X\otimes I}\uparrow \\
E\otimes\omega(X)\otimes\omega(Y) &                                          & E\otimes\omega(X)\otimes Y \\
\downarrow {\footnotesize I\otimes c_{\omega(X),\omega(Y)}} & & {\footnotesize c_{Y,\omega(X)}\otimes I}\uparrow \\
E\otimes\omega(Y)\otimes\omega(X) & \stackrel{e_Y\otimes I}{\longrightarrow} & E\otimes Y\otimes\omega(X)  
\end{array}\end{equation} 
Now we establish the relation between integral and idempotent objects.
\begin{lem}\label{intid}
For any integral object $(E,e)$ the pair $(E,e_E)$ is an idempotent object with $M = \omega(E)$ and $\Phi = \omega(e_E)$. 
\end{lem}
Proof:

The modified pentagon equation can be glued from two diagrams: 
$$\begin{array}{ccc}
E\otimes\omega(E)\otimes\omega(E) & \stackrel{e_E\otimes\omega(E)}{\longrightarrow} & E\otimes E\otimes\omega(E) \\
\uparrow {\footnotesize I\otimes\omega_{E,\omega(E)}} &                                       & || \\
E\otimes\omega(E\otimes\omega(E)) & \stackrel{e_{E\otimes\omega(E)}}{\longrightarrow} & E\otimes E\otimes\omega(E) \\
\downarrow {\footnotesize I\otimes\omega(e_E)} &                          & {\footnotesize I\otimes e_E}\downarrow \\
E\otimes\omega(E\otimes E)        & \stackrel{e_{E\otimes E}}{\longrightarrow} & E\otimes E\otimes E \\
\downarrow {\footnotesize I\otimes\omega_{E,E}} &                           & {\footnotesize e_E\otimes I}\uparrow \\
E\otimes\omega(E)\otimes\omega(E) &                                          & E\otimes\omega(E)\otimes E \\
\downarrow {\footnotesize I\otimes c_{\omega(E),\omega(E)}} &    & {\footnotesize c_{E,\omega(E)}\otimes I}\uparrow \\
E\otimes\omega(E)\otimes\omega(E) & \stackrel{e_E\otimes I}{\longrightarrow} & E\otimes E\otimes\omega(E)  
\end{array}$$ 
One (bottom part) is the diagram \ref{dinto} for $X=Y=E$ and the other (top part) is a functoriality diagram for collection $e$.
To complete the proof it is enough to note that composition of three top left vertical arrows in the glued diagram coincides with $I\otimes\Phi$. $\Box$ 

In particular, any integral object $(E,e)$ in the category $\G$ defines solution to pentagon equation $\Phi = \omega(e_E)$ and monoidal functor $\G_\Phi\to\G$ sending $k$ to $E$. 

Now we describe integral objects in categories of modules over Hopf algebras. 
Denote by $\M (H)$ the category of modules over Hopf algebra $H$. The forgetful functor $\M (H)\to\V{\it ect}_k$ has natural monoidal structure. 
\begin{prop}
Integral object in the category of modules over Hopf algebra with respect to forgetful functor are Hopf modules. 
\end{prop}
Sketch of the proof:

Let $V$ be a Hopf module. It can be checked directly that for any (left) $H$-module $X$ the map $V\otimes\omega(X)\to V\otimes X$ defined by $v\otimes x\mapsto \sum_{(v)}v_{(0)}\otimes v_{(1)}x$ is a structure of integral object.

Conversely, if $(E,e)$ is an integral object in $\M (H)$ then the composition 
$$E\stackrel{I\otimes i}{\longrightarrow} E\otimes\omega(H)\stackrel{e_H}{\longrightarrow} E\otimes H$$
is an $H$-linear comodule structure (Hopf module structure).
Here we consider Hopf algebra $H$ as a module over itself; $i:k\to\omega(H)$ denotes the unital inclusion $c\mapsto c1$. $\Box$

Combining the lemma and proposition above we can see that any Hopf module $V$ defines a solution to pentagon equation $\Phi$ on $M=\omega(V)$ (which is of course one defined in section \ref{esmpe}) and a monoidal functor $\G_\Phi\to\M(H)$ sending $k$ to $V$. The categorical reformulation of theorem \ref{spe} says that for any solution $\Phi$ of pentagon equation there is a monoidal functor from $\G_\Phi$ to the category of modules $\M(H)$ over some Hopf algebra $H$ sending $k$ to some Hopf module $V$. We outline the direct proof of this fact. First of all we need to define the category $\M(H)$ in terms of $\Phi$. 

Following \cite{sask} we call the pair $(X,\Psi)$ of a vector space $X$ and an isomorphism 
$v:M\otimes X\to M\otimes X$ a {\em module} over the solution to pentagon equation $(M,\Phi)$ if the following equality for operators on $M\otimes M\otimes X$ holds: 
\begin{equation}\label{se}
\Phi_{12}\Psi_{13}\Psi_{23} = \Psi_{23}\Phi_{12}.
\end{equation}
For example, the pair $(M,\Phi)$ is a module over itself and the pair $(X,I_{M\otimes X})$ is a module 
over $(M,\Phi)$ for any vector space $X$.  

Note that if we are given by Hopf $H$-module structure on $M$ defining $\Phi$ then any $H$-module becomes an $(M,\Phi)$-module by the formula $\Psi (m\otimes x) = \sum_{(m)}m_{(0)}\otimes m_{(1)}x$. 

A {\em morphism} of $(M,\Phi)$-modules $(X,\Psi)\to(X',\Psi')$ is a linear map $f:X\to X'$ such that 
\begin{equation}\label{he}
\Psi'(I\otimes f) = (I\otimes f)\Psi.
\end{equation}
It was shown in \cite{st} (in much more general setting) that the category $\M (M,\Phi)$ of modules over $(M,\Phi)$ forms rigid 
monoidal category with tensor product defined by
$$(X,\Psi)\otimes (X',\Psi') = (X\otimes X',\Psi_{12}{\Psi'}_{13}).$$
Forgetting of module structure gives a $k$-linear monoidal functor $\M (M,\Phi)\to\V{\it ect}_k$. 
Tannaka-Krein theory \cite{dm} decomposes this forgetful functor 
into an equivalence $\M (M,\Phi )\to \M (H)$ and a forgetful functor $\M (H)\to \V{\it ect}_k$ for category $\M (H)$ of (finite-dimensional over $k$) modules over some Hopf algebra $H$. It is not hard to check that this Hopf algebra coincides with one constructed in section \ref{sspe}. 

Now we concentrate on the role of Hopf bimodules in the picture. We prove that $(M,\Phi)$ has a natural structure of integral object in the category $\M (M,\Phi)$. 
\begin{prop}\label{iorc}
For any $(M,\Phi)$-module $(X,\Psi)$ structural map $\Psi$ defines the isomorphism of $(M,\Phi)$-modules
$$\Psi:(M,\Phi)\otimes X\to (M,\Phi)\otimes (X,\Psi).$$ 
In other words we have a functorial (in $(X,\Psi)$) collection of isomorphisms $e_{(X,\Psi)}:\Psi:(M,\Phi)\otimes\omega(X,\Psi)\to (M,\Phi)\otimes(X,\Psi).$

Moreover the diagrams of the form \ref{dinto} are commutative for this collection. 
\end{prop}
Sketch of the proof:

Indeed, the structural maps of the tensor products $(M,\Phi)\otimes X$, 
$(M,\Phi)\otimes (X,\Psi)$ are $\Phi_{12}$ and 
$\Phi_{12}\Psi_{13}$ correspondently. Hence we can consider the equation (\ref{se}) 
as a homomorphism condition (\ref{he}) for $\Psi$.

The commutativity of the diagram \ref{dinto} is equivalent to the decomposition of the equation 
\begin{equation}\label{setp}
\Phi_{12}\Psi_{13}{\Psi'}_{14}\Psi_{23}{\Psi'}_{24} = \Psi_{23}{\Psi'}_{24}\Phi_{12}
\end{equation}
into commutativity equation and equations of the form \ref{se} (up to some common factors):

using ${\Psi'}_{14}\Psi_{23} = \Psi_{23}{\Psi'}_{14}$ we replace left side of (\ref{setp}) by
$\Phi_{12}\Psi_{13}\Psi_{23}{\Psi'}_{14}{\Psi'}_{24}$,
then using $\Phi_{12}\Psi_{13}\Psi_{23} = \Psi_{23}\Phi_{12}$ we rewrite the result as
$\Psi_{23}\Phi_{12}{\Psi'}_{14}{\Psi'}_{24}$ and finally using $\Phi_{12}{\Psi'}_{13}{\Psi'}_{23} = {\Psi'}_{23}\Phi_{12}$ we get the right side of (\ref{setp}).
$\Box$

Combining lemma \ref{intid} with proposition \ref{iorc} we get a monoidal functor from the category $\G_\Phi$ to the category of modules $\M (M,\Phi )$ sending $k$ to $(M,\Phi )$. Tannaka-Krein equivalence $\M (M,\Phi )\cong \M(H)$ identify integral object $(M,\Phi )$ with some Hopf $H$-module. 

We finish this section reformulating connection between Hopf-Galois algebras and solutions to modified pentagon equation (and thus by proposition \ref{mofu} monoidal functors from category $\G_\Phi$ to vector spaces). It was proved in \cite{ul} that any linear monoidal functor from category $\M (H)$ of (finite-dimensional over $k$) modules over some Hopf algebra $H$ to vector spaces corresponds to some Hopf Galois $H$-module coalgebra. 
Thus the connection between Galois coalgebras and monoidal functors from $\G_\Phi$ to $\V{\it ect}_k$ can be represented by commutative diagram of monoidal functors: 
$$\begin{array}{ccc}
\G (\Phi ) & \stackrel{\omega_F}{\longrightarrow} & \V{\it ect}_k \\
\downarrow &                 & || \\
\M (H) & \stackrel{\omega_L}{\longrightarrow} & \V{\it ect}_k \\
\end{array}$$
where $\omega_L$ is a monoidal functor corresponding to Galois coalgebra $L$ and $\omega_F$ is a monoidal functor corresponding to solution $F$ of modified pentagon equation which is defined by $L$ (see example \ref{galmpe}).

\end{document}